\documentclass{amsart}

\advance\textheight by \topskip
\newtheorem{theorem}{Theorem}[section]
\newtheorem{corollary}[theorem]{Corollary}

\numberwithin{equation}{section}

\allowbreak

\title{Apostol-Euler polynomials arising from umbral calculus}

\author{Taekyun Kim}
\address{Department of Mathematics, Kwangwoon University, Seoul, S. Korea}
\email{tkkim@kw.ac.kr}

\author{Toufik Mansour}
\address{Department of Mathematics, University of Haifa, 31905 Haifa, Israel}
\email{tmansour@univ.haifa.ac.il}

\author{Seog-Hoon Rim}
\address{Department of Mathematics Education, Kyungpook National University, Taegu, S. Korea}
\email{shrim@knu.ac.kr}

\author{Sang-Hun Lee}
\address{Division of General Education, Kwangwoon University, Seoul, S. Korea}
\email{leesh58@kw.ac.kr}


\begin{document}

\begin{abstract}
In this paper, by using the orthogonality type as defined in the umbral calculus, we derive explicit formula for several well known polynomials as a linear combination of the Apostol-Euler polynomials.
\end{abstract}
\maketitle

\section{Introduction}
Let $\Pi_n$ be the set of all polynomials in a single variable $x$ over the complex field $\mathbb{C}$ of degree at most $n$. Clearly, $\Pi_n$ is a $(n+1)$-dimensional vector space over $\mathbb{C}$. Define
\begin{align}
\mathcal{H}=\left\{f(t)=\sum_{k\geq0} a_k\frac{t^k}{k!}\mid a_k\in\mathbb{C}\right\}.\label{eq7}
\end{align}
to be the algebra of formal power series in a single variable $t$. As is known, $\langle L|p(x)\rangle$ denotes the action of a linear functional $L\in \mathcal{H}$ on a polynomial $p(x)$ and we remind that the vector space on $\Pi_n$ are defined by
$$\langle cL+c'L'|p(x)\rangle=c\langle L|p(x)\rangle+c'\langle L'|p(x)\rangle,$$
for any $c,c'\in\mathbb{C}$ and $L,L'\in \mathcal{H}$ (see \cite{K4,K5,Ro1,Ro2}). The formal power series in variable $t$ define a linear functional on $\Pi_n$ by setting
\begin{align}
\langle f(t)|x^n\rangle=a_n, \mbox{ for all $n\geq0$, (see \cite{K4,K5,Ro1,Ro2})}.\label{eq8}
\end{align}
By (\ref{eq7}) and (\ref{eq8}), we have
\begin{align}
\langle t^k|x^n\rangle=n!\delta_{n,k}, \mbox{ for all $n,k\geq0$, (see \cite{K4,K5,Ro1,Ro2})},\label{eq9}
\end{align}
where $\delta_{n,k}$ is the Kronecker's symbol. Let $f_L(t)=\sum_{k\geq0}\langle L|x^k\rangle \frac{t^k}{k!}$ with $L\in \mathcal{H}$. From (\ref{eq9}), we have $\langle f_L(t)|x^n\rangle=\langle L|x^n\rangle$. So, the map $L\mapsto f_L(t)$ is a vector space isomorphic from $\Pi_n$ onto $\mathcal{H}$. Henceforth, $\mathcal{H}$ is thought of as set of both formal power series and linear functionals. We call $\mathcal{H}$ the {\em umbral algebra}. The {\em umbral calculus} is the study of umbral algebra.

Let $f(t)\in\mathcal{H}$.  The smallest integer $k$ for which the coefficient of $t^k$ does not vanish is called the {\em order} of $f(t)$ and is denoted by $O(f(t))$ (see \cite{K4,K5,Ro1,Ro2}). If $O(f(t))=1$, $O(f(t))=0$, then $f(t)$ is called a {\em delta}, an {\em invertable} series, respectively. For given two power series $f(t),g(t)\in\mathcal{H}$ such that $O(f(t))=1$ and $O(g(t))=0$, there exists a unique sequence $S_n(x)$ of polynomials with $\langle g(t)(f(t))^k|S_n(x)\rangle=n!\delta_{n,k}$ (this condition sometimes is called {\em orthogonality type}) for all $n,k\geq0$. The sequence $S_n(x)$ is called the {\em Sheffer} sequence for $(g(t),f(t))$ which is denoted by $S_n(x)\sim(g(t),f(t))$ (see \cite{K4,K5,Ro1,Ro2}).

For $f(t)\in\mathcal{H}$ and $p(x)\in\Pi$, we have
\begin{align}
\langle e^{yt}|p(x)\rangle=p(y),\quad \langle f(t)g(t)|p(x)\rangle=\langle f(t)|g(t)p(x)\rangle,\label{eq10}
\end{align}
and
\begin{align}
f(t)=\sum_{k\geq0}\langle f(t)|x^k\rangle\frac{t^k}{k!},\quad p(x)=\sum_{k\geq0}\langle t^k|p(x)\rangle\frac{x^k}{k!},\label{eq11}
\end{align}
(see \cite{K4,K5,Ro1,Ro2}). From (\ref{eq11}), we derive
\begin{align}
\langle t^k|p(x)\rangle=p^{(k)}(0), \langle 1|p^{(k)}(x)\rangle=p^{(k)}(0),\label{eq12}
\end{align}
where $p^{(k)}(0)$ denotes the $k$-th derivative of $p(x)$ respect to $x$ at $x=0$. Let $S_n(x)\sim(g(t),f(t))$. Then we have
\begin{align}
\frac{1}{g(\bar{f}(t))}e^{y\bar{f}(t)}=\sum_{k\geq0}S_k(y)\frac{t^k}{k!}, \label{eq14}
\end{align}
for all $y\in\mathbb{C}$, where $\bar{f}(t)$ is the compositional inverse of $f(t)$ (see \cite{K4,K5,K,R,Ro1,Ro2}).

For $\lambda\in\mathbb{C}$ with $\lambda\neq-1$, the {\em Apostol-Euler polynomials} (see \cite{H,N,T,W}) are defined by the generating function to be
\begin{align}
\frac{2}{\lambda e^t+1}e^{xt}=\sum_{k\geq0}E_k(x|\lambda)\frac{t^k}{k!}. \label{eqa1}
\end{align}
In particular, $x=0$, $E_n(0|\lambda)=E_n(\lambda)$ is called the {\em $n$-th Apostol-Euler number}. From (\ref{eqa1}), we can derive
\begin{align}
E_n(x|\lambda)=\sum_{k=0}^n\binom{n}{k}E_{n-k}(\lambda)x^k. \label{eqa2}
\end{align}
By (\ref{eqa2}), we have $\frac{d}{dx}E_n(x|\lambda)=nE_{n-1}(x|\lambda)$. Also, from (\ref{eqa1}) we have
\begin{align}
\frac{2}{\lambda e^t+1}=e^{E(\lambda)t}=\sum_{n\geq0}E_n(\lambda)\frac{t^n}{n!} \label{eqa3}
\end{align}
with the usual convention about replacing $E^n(\lambda)$ by $E_n(\lambda)$. By (\ref{eqa3}), we get
$$2=e^{E(\lambda)t}(\lambda e^t+1)=\lambda e^{(E(\lambda)+1)t}+e^{E(\lambda)t}=\sum_{n\geq0}(\lambda(E(\lambda)+1)^n+E_n(\lambda))\frac{t^n}{n!}.$$
Thus, by comparing the coefficients of the both sides, we have
\begin{align}
\lambda(E(\lambda)+1)^n+E_n(\lambda)=2\delta_{n,0}. \label{eqa4}
\end{align}
As is well known, {\em Bernoulli polynomial} (see \cite{AA,BK1,C,DD}) is also defined by the generating function to be
\begin{align}
\frac{t}{e^t-1}e^{xt}=\sum_{k\geq0}B_k(x)\frac{t^k}{k!}. \label{eqa5}
\end{align}

In the special case, $x=0$, $B_n(0)=B_n$ is called the {\em $n$-th Bernoulli number}. By (\ref{eqa5}), we get
\begin{align}
B_n(x)=\sum_{k=0}^n\binom{n}{k}B_{n-k}x^k. \label{eqa6}
\end{align}
From (\ref{eqa5}), we note that
\begin{align}
\frac{t}{e^t-1}=e^{Bt}=\sum_{n\geq0}B_n\frac{t^n}{n!} \label{eqa7}
\end{align}
with the usual convention about replacing $B^n$ by $B_n$. By (\ref{eqa6}) and (\ref{eqa7}), we get
$$t=e^{Bt}(e^t-1)=e^{(B+1)t}-e^{Bt}=\sum_{n\geq0}((B+1)^n-B_n)\frac{t^n}{n!},$$
which implies
\begin{align}
B_n(1)-B_n=(B+1)^n-B_n=\delta_{n,1},\quad B_0=1. \label{eqa8}
\end{align}
{\em Euler polynomials} (see \cite{AA,BK2,C,Ro2}) are defined by
\begin{align}
\frac{t}{e^t+1}e^{xt}=\sum_{k\geq0}E_k(x)\frac{t^k}{k!}. \label{eqa9}
\end{align}
In the special case, $x=0$, $E_n(0)=E_n$ is called the {\em $n$-th Euler number}. By (\ref{eqa9}), we get
\begin{align}
\frac{2}{e^t+1}=e^{Et}=\sum_{n\geq0}E_n\frac{t^n}{n!} \label{eqa10}
\end{align}
with the usual convention about replacing $E^n$ by $E_n$. By (\ref{eqa9}) and (\ref{eqa10}), we get
$$2=e^{Et}(e^t+1)=e^{(E+1)t}+e^{Et}=\sum_{n\geq0}((E+1)^n+E_n)\frac{t^n}{n!},$$
which implies
\begin{align}
E_n(1)+E_n=(E+1)^n+E_n=2\delta_{n,0}. \label{eqa11}
\end{align}
For $\lambda\in\mathbb{C}$ with $\lambda\neq-1$, the {\em Frobenius-Euler} (see \cite{AA,Cb,Ca, CCKS,CKOS}) polynomials are defined by
\begin{align}
\frac{1+\lambda}{e^t+\lambda}e^{xt}=\sum_{k\geq0}F_k(x|-\lambda)\frac{t^k}{k!}. \label{eqa12}
\end{align}
In the special case, $x=0$, $F_n(0|-\lambda)=F_n(-\lambda)$ is called the {\em $n$-th Frobenius-Euler number} (see \cite{Ca}). By (\ref{eqa12}), we get
\begin{align}
\frac{1+\lambda}{e^t+\lambda}=e^{Ft}=\sum_{n\geq0}F_n(-\lambda)\frac{t^n}{n!} \label{eqa13}
\end{align}
with the usual convention about replacing $F^n(-\lambda)$ by $F_n(-\lambda)$ (see \cite{Ca}). By (\ref{eqa12}) and (\ref{eqa13}), we get
$$1+\lambda=e^{F(-\lambda)t}(e^t+\lambda)=e^{(F(-\lambda)+1)t}+\lambda e^{F(-\lambda)t}=\sum_{n\geq0}((F(-\lambda)+1)^n+\lambda F_n(-\lambda))\frac{t^n}{n!},$$
which implies
\begin{align}
\lambda F_n(-\lambda)+F_n(1|-\lambda)=\lambda F_n(-\lambda)+(F(-\lambda)+1)^n=(1+\lambda)\delta_{n,0}. \label{eqa14}
\end{align}

In the nest section we present our main theorem and its applications. More precisely, by using the orthogonality type, we write any polynomial in $\Pi_n$ as a linear combination of the Apostol-Euler polynomials. Several applications related to Bernoulli, Euler and Frobenius-Euler polynomials are derived.

\section{Main results and applications}
Note that the set of the polynomials $E_0(x|\lambda),E_n(x|\lambda),\ldots,E_n(x|\lambda)$ is a good basis for $\Pi_n$. Thus, for $p(x)\in\Pi_n$, there exist constants $c_0,c_1,\ldots,c_n$ such that $p(x)=\sum_{k=0}^nc_kE_k(x|\lambda)$. Since $E_n(x|\lambda)\sim((1+\lambda e^t)/2,t)$ (see (\ref{eq14}) and \eqref{eqa1})), we have
$$\left\langle\frac{1+\lambda e^t}{2}t^k|E_n(x|\lambda)\right\rangle=n!\delta_{n,k},$$
which gives
$$\left\langle\frac{1+\lambda e^t}{2}t^k|p(x)\right\rangle=\sum_{\ell=0}^nc_\ell\left\langle\frac{1+\lambda e^t}{2}t^k|E_\ell(x|\lambda)\right\rangle=\sum_{\ell=0}^nc_\ell \ell!\delta_{\ell,k}=k!c_k.$$
Hence, we can state the following result.

\begin{theorem}\label{th1}
For all $p(x)\in\Pi_n$, there exist constants $c_0,c_1,\ldots,c_n$ such that $p(x)=\sum_{k=0}^nc_kE_k(x|\lambda)$, where
$$c_k=\frac{1}{2k!}\langle (1+\lambda e^t)t^k|p(x)\rangle.$$
\end{theorem}

Now, we present several applications for our theorem. As a first application, let us take $p(x)=x^n$ with $n\geq0$. By Theorem \ref{th1}, we have
$x^n=\sum_{k=0}^nc_kE_k(x|\lambda)$, where
\begin{align*}
c_k&=\frac{1}{2k!}\langle(1+\lambda e^t)t^k|x^n\rangle=\frac{1}{2}\binom{n}{k}\langle 1+\lambda e^t|x^{n-k}\rangle\\
&=\frac{1}{2}\binom{n}{k}(\delta_{n-k,0}+\lambda),
\end{align*}
which implies the following identity.

\begin{corollary}
For all $n\geq0$,
$$x^n=\frac{1}{2}E_n(x|\lambda)+\frac{\lambda}{2}\sum_{k=0}^n\binom{n}{k}E_k(x|\lambda).$$
\end{corollary}

Let $p(x)=B_n(x)\in\Pi_n$, then by Theorem \ref{th1} we have that
$B_n(x)=\sum_{k=0}^nc_kE_k(x|\lambda)$, where
\begin{align*}
c_k&=\frac{1}{2k!}\langle(1+\lambda e^t)t^k|B_n(x)\rangle=\frac{1}{2}\binom{n}{k}\langle 1+\lambda e^t|B_{n-k}(x)\rangle\\
&=\frac{1}{2}\binom{n}{k}(B_{n-k}+\lambda B_{n-k}(1)),
\end{align*}
which, by (\ref{eqa8}), implies the following identity.

\begin{corollary}
For all $n\geq2$,
$$B_n(x)=\frac{(\lambda-1)n}{4}E_{n-1}(x|\lambda)+\frac{1+\lambda}{2}\sum_{k=0,k\neq n-1}^n\binom{n}{k}B_{n-k}E_k(x|\lambda).$$
\end{corollary}

Let $p(x)=E_n(x)$, then by Theorem \ref{th1} we have that
$E_n(x)=\sum_{k=0}^nc_kE_k(x|\lambda)$, where
\begin{align*}
c_k&=\frac{1}{2k!}\langle(1+\lambda e^t)t^k|E_n(x)\rangle=\frac{1}{2}\binom{n}{k}\langle 1+\lambda e^t|E_{n-k}(x)\rangle\\
&=\frac{1}{2}\binom{n}{k}(E_{n-k}+\lambda E_{n-k}(1)),
\end{align*}
which, by (\ref{eqa11}), implies the following identity.

\begin{corollary}
For all $n\geq0$,
$$E_n(x)=\frac{1+\lambda}{2}\sum_{k=0}^n\binom{n}{k}E_{n-k}E_k(x|\lambda).$$
\end{corollary}

For another application, let $p(x)=F_n(x|-\lambda)$, then by Theorem \ref{th1} we have that
$F_n(x|-\lambda)=\sum_{k=0}^nc_kE_k(x|\lambda)$, where
\begin{align*}
c_k&=\frac{1}{2k!}\langle(1+\lambda e^t)t^k|F_n(x|-\lambda)\rangle=\frac{1}{2}\binom{n}{k}\langle 1+\lambda e^t|F_{n-k}(x|-\lambda)\rangle\\
&=\frac{1}{2}\binom{n}{k}(F_{n-k}(-\lambda)+\lambda F_{n-k}(1|-\lambda)).
\end{align*}
which, by (\ref{eqa14}), implies the following identity.

\begin{corollary}
For all $n\geq1$,
$$F_n(x|-\lambda)=\frac{1+\lambda}{2}E_n(x|\lambda)+\frac{1-\lambda^2}{2}\sum_{k=0}^{n-1}\binom{n}{k}F_{n-k}(-\lambda)E_k(x|\lambda).$$
\end{corollary}

Again, let $p(x)=y_n(x)=\sum_{k=0}^n \frac{(n+k)!}{(n-k)!k!}\frac{x^k}{2^k}$ be the $n$-th {\em Bessel polynomial} (which is the solution of the following differential equation $x^2f''(x)+2(x+1)f'+n(n+1)f=0$, where $f'(x)$ denotes the derivative of $f(x)$, see \cite{Ro1,Ro2}). Then by Theorem \ref{th1}, we can write  $y_n(x)=\sum_{k=0}^nc_kE_k(x|\lambda)$, where
\begin{align*}
c_k&=\frac{1}{2k!}\sum_{\ell=0}^n \frac{(n+\ell)!}{(n-\ell)!\ell!2^\ell} \langle 1+\lambda e^t|t^kx^\ell\rangle\\
&=\frac{1}{2}\sum_{\ell=k}^n \frac{(n+\ell)!}{(n-\ell)!\ell!2^\ell}\binom{\ell}{k} \langle 1+\lambda e^t|x^{\ell-k}\rangle\\
&=\frac{1}{2}\sum_{\ell=k}^n \frac{(n+\ell)!}{(n-\ell)!\ell!2^\ell}\binom{\ell}{k} (\delta_{n-k,0}+\lambda)\\
&=\frac{k!}{2^{k+1}}\binom{n}{k}\binom{n+k}{k}+
\lambda\sum_{\ell=k}^n \frac{k!}{2^{\ell+1}}\binom{\ell}{k}\binom{n}{\ell}\binom{n+\ell}{\ell}\\
\end{align*}
which implies the following identity.

\begin{corollary}
For all $n\geq1$,
$$y_n(x)=\sum_{k=0}^{n}\frac{k!}{2^{k+1}}\binom{n}{k}\binom{n+k}{k}E_k(x|\lambda)+\lambda\sum_{k=0}^{n}\sum_{\ell=k}^n \frac{k!}{2^{\ell+1}}\binom{\ell}{k}\binom{n}{\ell}\binom{n+\ell}{\ell}E_k(x|\lambda).$$
\end{corollary}

We end by noting that if we substitute $\lambda=0$ in any of our corollaries, then we get the well known value of the polynomial $p(x)$. For instance, by setting $\lambda=0$, the last corollary gives that $y_n(x)=\sum_{k=0}^n \frac{(n+k)!}{(n-k)!k!}\frac{x^k}{2^k}$ as expected.


\end{document}